# DISCUSSION: A TALE OF THREE COUSINS: LASSO, L2BOOSTING AND DANTZIG

By N. Meinshausen,[1] G. Rocha[2] and B. Yu[3]

*University of California, Berkeley*

We would like to congratulate the authors for their thought-provoking and interesting paper. The Dantzig paper is on the timely topic of high-dimensional data modeling that has been the center of much research lately and where many exciting results have been obtained. It also falls in the very hot area at the interface of statistics and optimization: $\ell_1$-constrained minimization in linear models for computationally efficient model selection, or sparse model estimation (Chen, Donoho and Saunders [5] and Tibshirani [17]). The sparsity consideration indicates a trend in high-dimensional data modeling advancing from prediction, the hallmark of machine learning, to sparsity—a proxy for interpretability. This trend has been greatly fueled by the participation of statisticians in machine learning research. In particular, Lasso (Tibshirani [17]) is the focus of many sparsity studies in terms of both theoretical analysis (Knight and Fu [10], Greenshtein and Ritov [9], van de Geer [19], Bunea, Tsybakov and Wegkamp [3], Meinshausen and Bühlmann [13], Zhao and Yu [23] and Wainwright [20]) and fast algorithm development (Osborne, Presnell and Turlach [15] and Efron et al. [8]).

Given $n$ units of data $Z_i = (X_i, Y_i)$ with $Y_i \in \mathbb{R}$ and $X_i^T \in \mathbb{R}^p$ for $i = 1, \ldots, n$, let $Y = (Y_1, \ldots, Y_n)^T \in \mathbb{R}^n$ be the continuous vector response variable and $X = (X_1, \ldots, X_n)^T$ the $n \times p$ design matrix and let the columns of $X$ be normalized to have $\ell_2$-norm 1. It is often useful to assume a linear regression model,

$$(1) \qquad\qquad Y = X\beta + \varepsilon,$$

where $\varepsilon$ is an i.i.d. $N(0, \sigma^2)$ vector of size $n$.

Received February 2007.

[1] Supported in part by a scholarship from DFG (Deutsche Forschungsgemeinschaft).

[2] Supported in part by ARO Grant W911NF-05-1-0104.

[3] Supported in part by a Guggenheim Fellowship, NSF Grants DMS-06-05165 and DMS-03-036508, and ARO Grant W911NF-05-1-0104.

Lasso minimizes the $\ell_1$-norm of the parameters subject to a constraint on squared error loss. That is, $\beta^{\text{lasso}}(t)$ solves the $\ell_1$-constrained minimization problem

$$(2) \qquad \min_{\beta} \|\beta\|_1 \quad \text{subject to} \quad \tfrac{1}{2}\|Y - X\beta\|_2^2 \leq t.$$

We can clearly use constraint and objective function interchangeably. For each value of $t > 0$, one can also find a value of the Lagrange multiplier $\lambda$ so that Lasso is the solution of the penalized version

$$(3) \qquad \min_{\beta} \tfrac{1}{2}\|Y - X\beta\|_2^2 + \lambda\|\beta\|_1.$$

Finally, it is well known that an alternative form of Lasso (Osborne, Presnell and Turlach [15]) asserts that $\beta_{\lambda}^{\text{lasso}}$ also solves

$$(4) \qquad \min_{\beta} \tfrac{1}{2}\beta^T X^T X\beta \quad \text{subject to} \quad \|X^T(Y - X\beta)\|_{\infty} \leq \lambda,$$

where $\lambda$ is identical to the penalty parameter in the penalized version (3). In what follows, we consider Dantzig estimates $\beta_{\lambda}^{\text{dantzig}}$ solving the constrained minimization problem

$$(5) \qquad \min_{\beta} \|\beta\|_1 \quad \text{subject to} \quad \|X^T(Y - X\beta)\|_{\infty} \leq \lambda.$$

The Dantzig selector as proposed by the authors uses $\lambda = \lambda_p(\sigma) = \sigma\sqrt{2\log p}$. To distinguish the two, we reserve the term Dantzig selector for this particular choice of $\lambda$ throughout this discussion. Comparing Dantzig with Lasso in its forms (4) and (5) reveals very clearly their close kinship. Hence we would like to view the Dantzig paper in the context of the vast literature on Lasso. We will start with some comments on the theoretical side before concentrating on comparing Dantzig and Lasso from the points of view of algorithmic and statistical performance.

**1. Lasso and Dantzig: theoretical results.** Assuming $\sigma$ is known, the Dantzig selector uses a fixed tuning parameter $\lambda_p(\sigma)$. Under a condition called *Uniform Uncertainty Principle* (which requires almost orthonormal predictors when choosing subsets of variables), an effective bound is obtained on the MSE $\|\hat{\beta}_{\lambda_p(\sigma)}^{\text{dantzig}} - \beta\|_2^2$ for the Dantzig selector. After a simple step of bounding the projected errors on the predictors, the proof is deterministic. This bounding step gives rise to the particular chosen threshold $\lambda_p(\sigma)$. In terms of tools used, this paper is closely related to earlier papers by the authors, Donoho, Elad and Temlyakov [7] and Donoho [6] on Lasso.

There is a parallel development of understanding Lasso under the linear regression model in (1) with stochastic tools. The results are in terms of the $\ell_2$-MSE on $\beta$ and also in terms of the $\ell_2$-MSE on the regression function $X\beta$ (e.g., Greenshtein and Ritov [9], Bunea, Tsybakov and Wegkamp [3], van de Geer [19], Zhang and Huang [22] and Meinshausen and Yu [14]). Related results for L2Boosting are obtained by Bühlmann [2]. Since Lasso is important for its model selection property, it is natural to study directly Lasso's model selection consistency as in the work of Meinshausen and Bühlmann [13], Zhao and Yu [23], Zou [24], Wainwright [20, 21] and Tropp [18]. What has emerged from this cluster of work is the necessity of an irrepresentable condition for Lasso to select the correct variables under sparsity conditions on the model. This condition regulates how correlated the predictors can be before wrong predictors are selected. However, this condition can be relaxed and Lasso still behaves sensibly. Specifically, Meinshausen and Yu [14] and Zhang and Huang [22] assume less restricted conditions on the predictors than the UUP condition to derive a bound on the same MSE ($\beta$) for an arbitrary $\lambda$. The bound is probably weaker than the Dantzig bound, but the assumptions are weaker as well so it covers commonly occurring highly correlated predictors. It is a consequence of this bound that in the case of $p \gg n$, if the model is sparse, Lasso can reduce significantly the number of predictors while keeping the correct ones. It would be interesting to see the Dantzig bound generalized to the case of more correlated predictors and for a range of $\lambda$'s since $\sigma$ is mostly unknown in practice and has to be estimated.

## 2. Lasso and Dantzig: algorithm and performance.
The similarities of Lasso and Dantzig revealed in (4) and (5) beg us to ask: How does Dantzig differ from Lasso? Which one should one use in practice and why? Let us start with a simple case where geometric visualizations of Dantzig and Lasso optimization problems can be easily displayed.

*Lasso versus Dantzig: $p = 3$ and in the population limit.* We choose three predictors from the multivariate normal distribution with a zero mean vector and a covariance matrix $V$ with a unit diagonal and entries $V_{12} = 0$ and $V_{13} = V_{23} = r$, where $|r| < 1/\sqrt{2}$ to guarantee positive definiteness of $V$. For simplicity, we consider the case of $n = \infty$, so we have zero noise and the population covariance $V$. We do this by setting the observations to be $Y = X\beta^*$, with $\beta^* = (1, 1, 0)$ and $X$ given by the Cholesky decomposition of $V$ so $X'X = V$. For the purpose of visualization, we rewrite the minimization problems in (2) and (5) in the alternative forms

$$(6) \qquad \min_{\beta} \|Y - X\beta\|_2^2 \quad \text{subject to} \quad \|\beta\|_1 \leq t, \quad \text{for Lasso;}$$

$$(7) \qquad \min_{\beta} \|X'(Y - X\beta)\|_\infty \quad \text{subject to} \quad \|\beta\|_1 \leq t, \quad \text{for Dantzig.}$$

In Figure 1, we display six plots of these alternative minimization problems. In the two leftmost columns, the $\ell_1$-polytopes sitting at the origin give the same $\ell_1$-constraint $\|\beta\|_1 \le 1$. The touching ball or ellipsoids in the first row correspond to the Lasso $\ell_2$-objective function for the Lasso, while the cube and polytopes in the second row correspond to the $\ell_\infty$-objective function for Dantzig. In the first column of the plots, $r = 0$ and both Lasso and Dantzig correctly select only the first two variables. In the second column, we set the correlation at $r = 0.5$. The Lasso still correctly selects only the first two variables. Meanwhile, the Dantzig admits multiple solutions, namely all points belonging to the line connecting $(0, 0, 1)$ and $\frac{(1,1,0)}{2}$. While it is true that $\frac{(1,1,0)}{2}$ is one of the Dantzig solutions correctly selecting the first two variables and discarding the third, all other solutions incorrectly include the third variable. In the other extreme, $(0, 0, 1)$ is also a solution where the first and second variables are wiped out from the model and only the third is added.

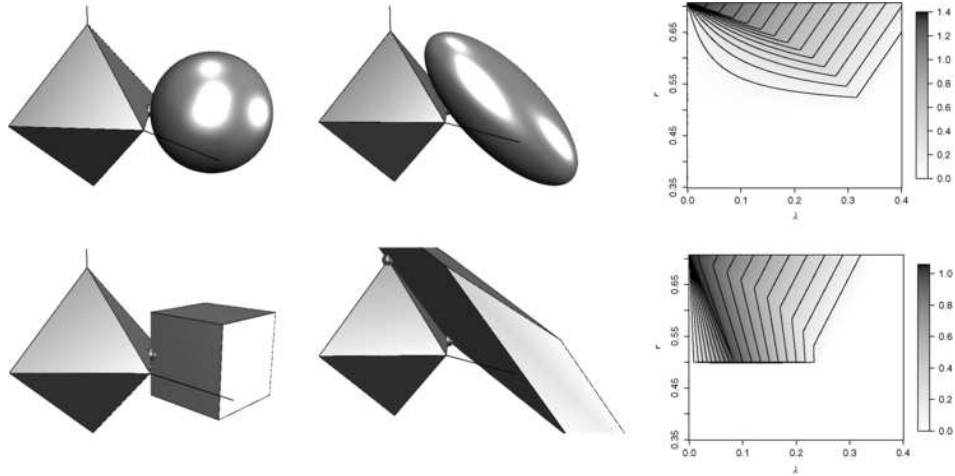

Fig. 1. *The panels in the first row and second row refer respectively to Lasso and Dantzig. The geometry in the $\beta$ space for the optimization problems (6) and (7) is shown for the uncorrelated design (leftmost panel) and for correlated design with $r = 0.5$ (middle panel). The Lasso solution is the point where the ellipsoid of $\ell_2$-loss touches the $\ell_1$-polytope and is unique in both cases. For Dantzig, the solutions are given by the points touching the $\ell_1$-polytope and the box-shaped $\ell_\infty$-constraint on the correlations of the predictor variables with the residuals. For $r = 0.5$, the solution is not uniquely determined for Dantzig as the side of the box aligns with the surface of the $\ell_1$-polytope. The rightmost column shows the third component $\hat{\beta}_3$ of the respective solution as a function of the correlation $r$ and the regularization parameter $\lambda$ as defined in (4) and (5). The Dantzig solution is not continuous at $r = 0.5$.*

In this example, $r = 0.5$ is a critical point where the irrepresentable condition (Zhao and Yu [23]) breaks down. The transition from below 0.5 to above can be seen in the third column of Figure 1, which depicts the contour plots of the estimated $\beta_3$ by Lasso and Dantzig: $r$ varies from 0.35 to 0.70 along the vertical direction and each horizontal line shows the whole path as a function of $\lambda$ in the optimization problems (4) and (5) for a fixed $r$. When $r > 0.5$, both Lasso and Dantzig systematically select the wrong third predictor (or the estimated $\beta_3$'s are nonzero). In terms of size of the incorrectly added coefficient, however, the transition is much sharper for Dantzig as $r$ crosses 0.5. In fact, the solution of the Dantzig is not a Lipschitz continuous function of the observations for $r = 0.5$. This could be expected, as Dantzig is the solution of a linear program (LP) problem and the estimator can thus jump from one vertex in the $\ell_\infty$ box to another if the data changes slightly. When $\lambda$ varies, the regularization path for the Dantzig is piecewise linear. However, the flat faces of both the loss and the penalty functions can cause jumps in the path, similarly to what happens in the $\ell_1$-penalized quantile regression (Rosset and Zhu [16]). This makes the design of an algorithm in the spirit of the homotopy/LARS-LASSO algorithm for the Lasso (Osborne, Presnell and Turlach [15] and Efron et al. [8]) more challenging and gives rise to jittery paths relative to Lasso and L2Boosting, as seen in the simulated example below.

The first column of Figure 1 suggests that Lasso and Dantzig could coincide. At the very least, their regularization paths share the same terminal points given by the minimal $\ell_1$-norm vector of coefficients, causing the correlation of all predictors with the residuals to be zero. In fact, more similarities exist: we now provide a sufficient condition for the two paths to entirely agree when $n \geq p$. The condition is diagonal dominance of $(X^T X)^{-1}$, that is, for $M = (X^T X)^{-1}$,

$$(8) \qquad M_{jj} > \sum_{i \neq j} |M_{ij}| \qquad \text{for all } j = 1, \ldots, p.$$

When $p = 2$, condition (8) is always satisfied so Lasso is exactly the same as Dantzig (and L2Boosting). Moreover, the irrepresentable condition is always satisfied as well. The diagonal dominance condition (8) is related to the positive cone condition used in Efron et al. [8] to show that L2Boosting and Lasso share the same path. The positive cone condition requires, for all subsets $\mathcal{A} \subseteq \{1, \ldots, p\}$ of variables, that $M_{jj} > -\sum_{i \neq j} M_{ij}$, where $M = (X_{\mathcal{A}}^T X_{\mathcal{A}})^{-1}$ and is always trivially satisfied for $p = 2$.

THEOREM 1. *Under the diagonal dominance condition (8), the Lasso solution (3) and the Dantzig solution (5) are identical for any value of $\lambda > 0$ (Lasso and Dantzig share the same path).*

PROOF. First, define the vector $g(\beta) = X^T(Y - X\beta) \in \mathbb{R}^p$ containing the correlation of the residuals with the original predictor variables. The Lasso solution is unique under condition (8). A necessary and sufficient condition for a vector $\beta$ to be the Lasso solution is, by the Karush–Kuhn–Tucker conditions (Bertsekas [1]), that (a) for all $k$: $g_k(\beta) \in [-\lambda, \lambda]$ and (b) for all $k \in \{l : \beta_l \neq 0\}$ it holds that $g_k(\beta) = \lambda \operatorname{sign}(\beta_k)$. We show that the Dantzig solution (5) is a valid Lasso solution under diagonal dominance (8). The Dantzig fulfills condition (a) by construction.

We now show that the (unique) Dantzig solution also satisfies (b). Assume to the contrary that $\beta$ is a solution of the Dantzig and there is some $j \in \{k : \beta_k \neq 0\}$ such that $g_j(\beta) \in [-\lambda, \lambda]$ but $g_j(\beta) \neq \lambda \operatorname{sign}(\beta_j)$. Let $\delta \in \mathbb{R}^p$ be a vector with $\delta_k = 0$ for all $k \neq j$ and $\delta_j = \operatorname{sign}(\beta_j)$ and define $\gamma = -(X^T X)^{-1} \delta$. We have $g(\beta + \nu\gamma) = g(\beta) + \nu\delta$, so only the $j$th component of the vector of correlations is changed by an amount $\nu \operatorname{sign}(\beta_j)$. Since we have assumed $|g_j(\beta)| < \lambda$, there exists some $\nu > 0$ such that $\beta + \nu\gamma$ is still in the feasible region.

To complete the proof we now show that, under the diagonal dominance condition (8), the $\ell_1$-norm of $\beta + \nu\gamma$ will be smaller than the $\ell_1$-norm of $\beta$ for small values of $\nu$. Denote by $\beta_{-j}$ the vector with entries identical to $\beta$, except for the $j$th component, which is set to zero. We can write

$$\|\beta + \nu\gamma\|_1 \leq \|\beta_{-j}\|_1 + \nu\|\gamma_{-j}\|_1 + |\beta_j + \nu\gamma_j|$$

$$\leq \|\beta_{-j}\|_1 + \nu \sum_{k \neq j} |M_{kj}| + |\beta_j| - \nu M_{jj}$$

$$= \|\beta\|_1 - \nu\left(M_{jj} - \sum_{k \neq j} |M_{kj}|\right),$$

where the first inequality results from using the triangle inequality twice and the second inequality stems from $\gamma_k = -M_{kj}\operatorname{sign}(\beta_j)$ with $M = (X^T X)^{-1}$. It thus holds that, for small enough values of $\nu > 0$, the right-hand side is smaller than $\|\beta\|_1$ under the diagonal dominance condition (8). Hence, the vector $\beta$ with $g_j(\beta) \neq \lambda \operatorname{sign}(\beta_j)$ cannot be the Dantzig solution. We conclude that the Dantzig solution must satisfy properties (a) and (b) and thus coincides with the Lasso solution (3). □

As alluded to earlier, the Dantzig selector needs the true $\sigma$ to be applied to real-world data. One obvious alternative is to use the Dantzig path and cross-validation. This gives another reason for obtaining the whole path. We define our *data-driven Dantzig selector* (*DD*) by computing $\hat{\sigma}_{\mathrm{CV}}^2$—the smallest fivefold cross-validated mean squared error over the Dantzig path— and plugging it into $\lambda_p(\hat{\sigma}_{\mathrm{CV}})$. Needless to say, this estimator is not without its problems: one being that the cross-validated error might not be a good estimate of the prediction error in the $p \gg n$ case and the other that it might overestimate $\sigma^2$. However, we decide to use it because it is sensible and simple. We later compare the performance of the data-driven Dantzig selector with the Dantzig estimator corresponding to the $\hat{\lambda}_{\mathrm{CV}}$ chosen as the minimizer of the cross-validated mean squared error.

A more realistic simulation example is in order for further comparisons of Lasso and Dantzig. The following simulation example reflects the common $p > n$ situation seen in recent real-world data applications. L2Boosting, Lasso and Dantzig will be contrasted against each other in terms of algorithmic and performance behavior. Path smoothness will be examined and statistical performance criteria include MSE on the $\beta$, MSE on the regression function $X\beta$ and a variable selection quality plot (i.e., correctly selected variables relative to falsely selected variables). In addition, we vary the signal to noise ratio and correlation level of the predictors to bring out more insight.

*Lasso, L2Boosting and Dantzig: $p > n$ and correlated predictors.* We consider random design with $p = 60$ variables and $n = 40$. Predictor variables have a multivariate Gaussian distribution $X \sim \mathcal{N}(0, \Sigma)$, where the population covariance matrix $\Sigma$ of the predictor variables is Toeplitz, that is, $\Sigma_{ij} = \rho^{|i-j|}$ for all $1 \leq i, j \leq p$. The response vector $Y$ is obtained as in (1),

$$Y = X\beta^* + \sigma\varepsilon, \tag{9}$$

where $\varepsilon = (\varepsilon_1, \ldots, \varepsilon_n)$ is i.i.d. noise with a standard Gaussian distribution. The $p$-dimensional vector $\beta^*$ is drawn once from a standard Gaussian distribution and all but 10 randomly selected coefficients are set to zero. To be precise, the true parameter vector $\beta^*$ used has entries

$$-0.65, -0.38, -0.37, -0.27, -0.12, -0.08, 0.05, 0.24, 0.37, 0.41,$$

for components 60, 2, 21, 49, 20, 27, 4, 43, 51, 32, with all other components set to zero. Three simulation setups are (a) $\rho = 0$, $\sigma = 0.2$; (b) $\rho = 0.9$, $\sigma = 0.2$; (c) $\rho = 0.9$, $\sigma = 0.6$. The vector $\beta^*$ is rescaled in each case so that $\|X\beta^*\|_2^2 = n$. We do not include the case that $\rho = 0$ and $\sigma = 0.6$ for the results are similar to (a).

Computing the solution path for both Lasso and L2Boosting took under half a second of CPU time each, using the LARS software in R of Efron et al. [8]. Computing the solution path of the Dantzig for 200 distinct values of the regularization parameter $\lambda$ took more than 30 seconds on the same computer, using either a standard C linear programming library lp_solve (called from R) or the Matlab code supplied in the $\ell_1$-magic package (Candès [4]). The relatively long running time for the current Dantzig algorithms makes it necessary to develop a path-following algorithm. As mentioned before, the Dantzig path could have jumps and, as a result, its path-following algorithm could be somewhat more involved, as in Li and Zhu [12].

Other simulations with different randomly chosen sparse $\beta^*$'s were conducted and yielded similar results as was demonstrated with this particular choice of $\beta^*$. In almost all cases, Lasso and L2Boosting outperform Dantzig

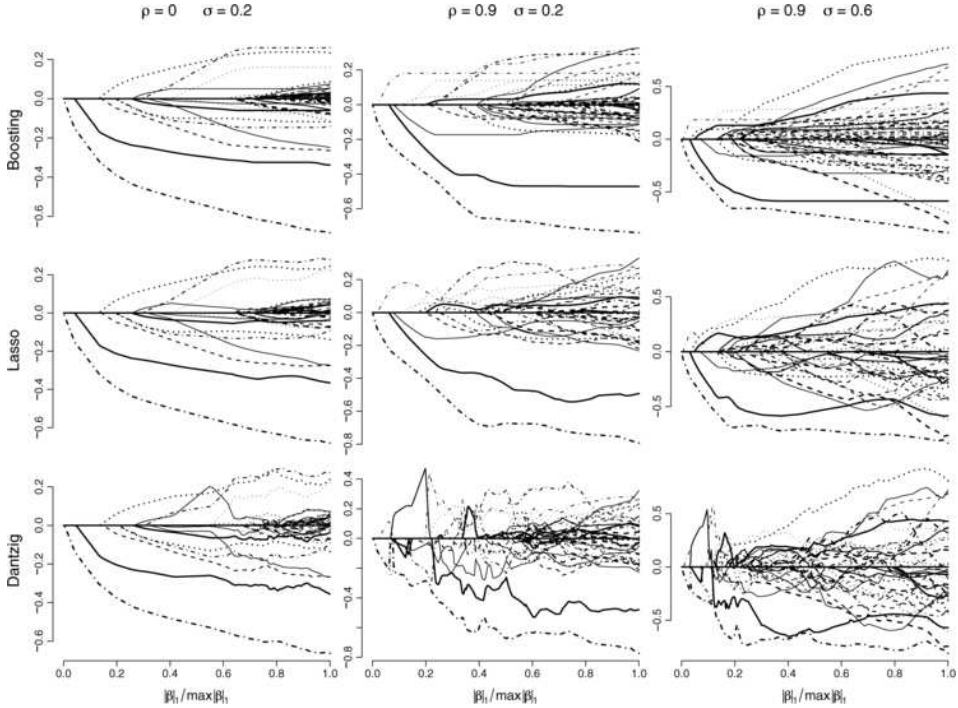

Fɪɢ. 2.   *Regularization paths from a single realization for each setup* (a), (b) *and* (c) *for L2Boosting (first row), Lasso (second row) and Dantzig (third row). The Dantzig path is jittery for a very correlated design (large value of* $\rho$*). The ends of the paths (for* $\lambda \to 0$*) agree for Dantzig and Lasso.*

and the Dantzig path is more jittery; when signal to noise ratio (SNR) is relatively high and the predictors are highly correlated, the performance gain of L2Boosting and Lasso over Dantzig cannot be ignored.

Now let us look into the details of the results in Figures 2, 3 and 4. Figure 2 displays path plots under (a), (b) and (c) for a single realization of the linear model (9). The horizontal axes are scaled so that the path plots are comparable. Given everything else being equal, a correlation increase or an SNR decrease makes the path more jittery for all three methods, with various degrees. Across methods, L2Boosting's path is most smooth, Lasso's is less smooth and Dantzig's is most jittery. Moreover, under the same simulation setup, the branching points from zero of the three methods are quite similar although the path smoothness differs.

Does the smoothness/jittery property of the path of a method readily translate into meaningful performance properties? Figures 3 and 4 attempt to answer this question. The first one shows that in terms of both MSE's, Lasso and L2Boosting are similar and in general better than Dantzig over the whole path. The improvement of Lasso or L2Boosting over Dantzig is more pronounced for the MSE on $\beta$ than that on $X\beta$. The middle column in Figure 3, with high correlation between predictors and high SNR, shows the worst case for Dantzig, relative to L2Boosting and Lasso. Such results are in terms of both MSE's, with the MSE for $\beta$ worse than the MSE for $X\beta$. This indicates qualitatively a regime where, when correlation and SNR are matched in some way, Dantzig is worse off than L2Boosting and Lasso. In other words, Lasso and L2Boosting are more effective to extract statistical information. With the same high correlation, however, when the SNR decreases (as shown in the right column of Figure 3), the statistical problem becomes hard for all of them and the advantage of Lasso and L2Boosting diminishes. For both MSE's, cross-validation selects better tuning parameters for all three methods than the data-driven Dantzig (DD) with the exception of setup (c). In this setup, the noise level is high and so is the correlation level, estimating individual $\beta$'s becomes difficult and hence it is better to be conservative as $\lambda_{DD}$ sets many $\beta$'s to zero (cf. the rightmost plot in the second row of Figure 3). However, when the performance measure is on prediction or the MSE on $X\beta$, $\lambda_{CV}$ does better again than $\lambda_{DD}$ (cf. the rightmost plot in the fourth row of Figure 3).

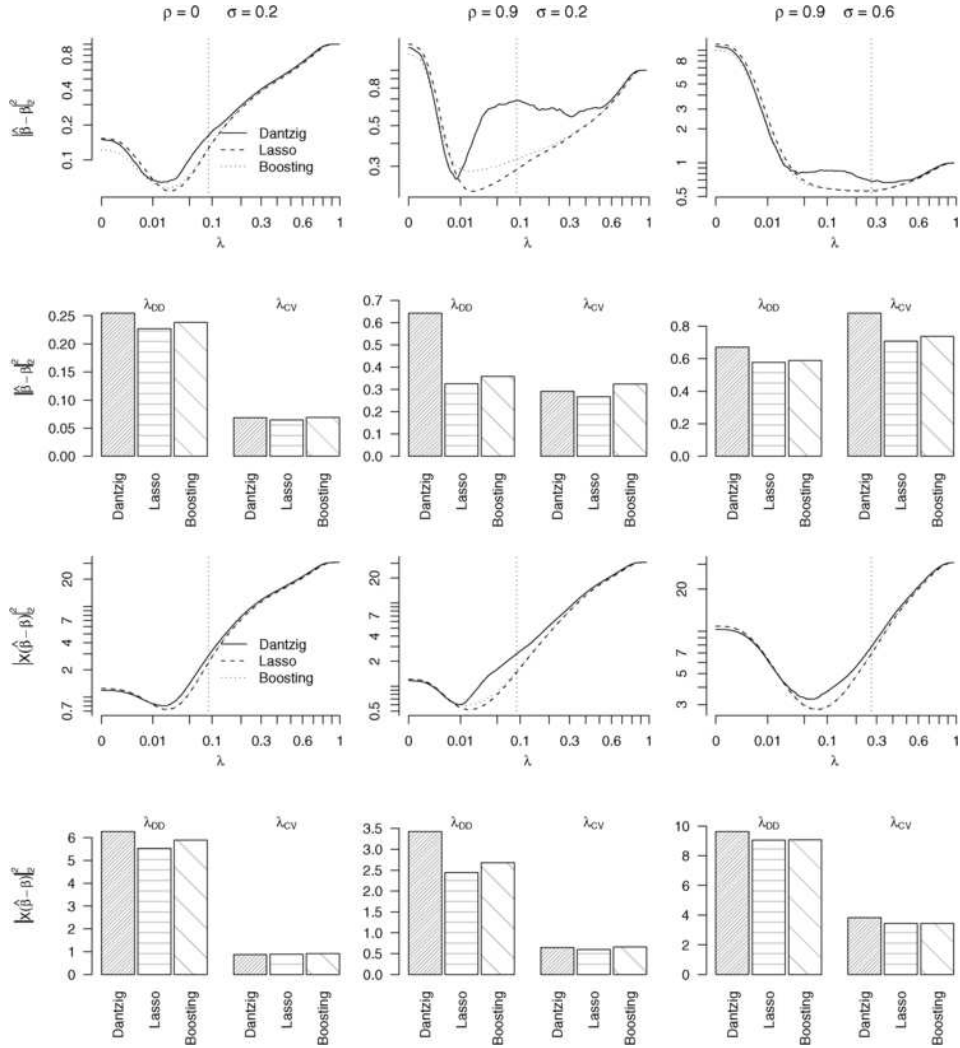

FIG. 3. *For the three setups* (a), (b) *and* (c), *the first row shows the MSE's on β of the Dantzig, the Lasso and L2Boosting solution as a function of the regularization parameter λ, averaged over 50 simulations. All three methods perform approximately equally well, with the exception of setting* (b), *where Dantzig performs worse. The vertical dotted line indicates the proposed fixed value of $λ_p(σ)$. The second row compares the solutions obtained by using the data-driven ($λ_{DD}$) and the cross-validation ($λ_{CV}$) tuning of the regularization parameter. In general, cross-validation gives a better fit except for the third setting* (c) *where the MSE on β favors the conservative data-driven Dantzig selector. The next two rows show comparable plots for the MSE's on Xβ. Here, the difference between all three methods is even smaller. For all three setups, the cross-validation tuned regularization parameter $λ_{CV}$ always results in a better MSE on Xβ or a better predictive performance than its data-driven Dantzig selector counterpart $λ_{DD}$.*

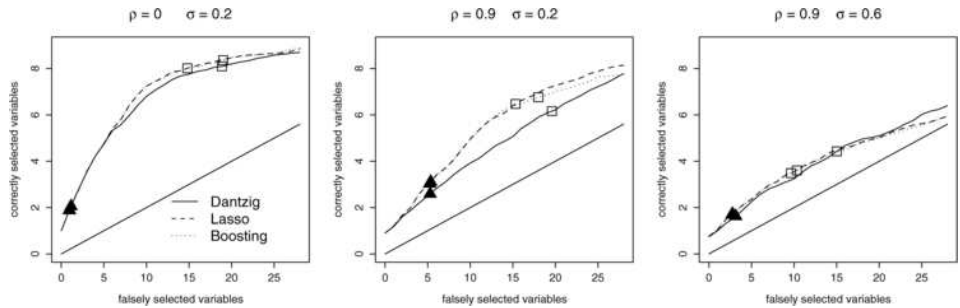

FIG. 4. *The average number of correctly selected variables as a function of the number of falsely selected variables, averaged over 50 simulations. The straight line corresponds to the performance under random selection of variables. Filled triangles indicate the solution under $\lambda_{DD}$, whereas the solution for $\lambda_{CV}$ is marked by squares.*

Last but not least, we assess the model selection prospect of the three methods with the CV-selected or the DD-selected tuning parameter $\lambda$. Figure 4 contains three plots under the three simulation setups. The horizontal axis plots the number of falsely selected variables and the vertical gives the corresponding correctly selected variables. Within each plot, the straight line gives the result of random selection of predictors; the solid curved line is Dantzig, dashed line is Lasso and dotted line is L2Boosting. The triangles indicate the DD selection and squares the CV selection of tuning parameters, for each method depending on the curve where the symbol is sitting. Obviously, all methods do better than random selection and the gain is highest when the predictors are not correlated. The gain is reduced when the correlation is high, but with a larger gain in the case of high SNR (middle plot) than the low-SNR case (right plot). In particular, the most differentiating case is setup (b): high correlation and high SNR. For all three methods, CV would pick up two or three more correct predictors with the same false predictors as random selection, and there is a slight but definite advantage of L2Boosting and Lasso over Dantzig. For high correlation and low SNR, only one or two correct ones can be gained over random selection of the same number of falsely selected predictors. Clearly, DD is very conservative to select very few predictors for all three methods, while CV has a tendency to include too many noise variables for low SNR; this is well known and has already been studied in more detail in Leng, Lin and Wahba [11] and Meinshausen and Bühlmann [13]. Nevertheless, for all three methods CV seems to give a better balance on the total number of correct predictors and false predictors. For any choice of the regularization parameter, L2Boosting and Lasso are in general no worse and sometimes better than Dantzig.

**3. Concluding remarks.** In this discussion, we have attempted to understand the Dantzig selector in relation to its cousins Lasso and L2Boosting. We believe that computing Dantzig or the Lasso for a single value of the penalty parameter $\lambda$ does not work well in practice; we need the entire solution path to select a meaningful model with good predictive performance. Without a path-following algorithm, computing the solution path for Dantzig is computationally very intensive (which is the reason we were limited to rather small data sets for the numerical examples). Leaving aside computational aspects, the first visual impression of the Dantzig solution path is its jitteriness when compared to the much smoother Lasso or L2Boosting solution paths, especially for highly correlated predictor variables. However, we showed that the smoothness of the path is not always indicative of performance. For the same regularization parameter, Lasso and L2Boosting performed in all settings at least as well as the Dantzig selector (and sometimes substantially better) and Dantzig performed on par with Lasso and L2Boosting for low signal to noise ratio even though its path is much more jittery. For almost all settings considered, the regularization parameter selected by cross-validation gives better MSE's than the data-driven Dantzig selector. In summary, we have not yet seen compelling evidence that would persuade us to use the Dantzig in practice rather than Lasso or L2Boosting.

**Acknowledgment.** We would like to thank Martin Wainwright for helpful comments on an earlier version of the discussion.

N. MEINSHAUSEN
DEPARTMENT OF STATISTICS
UNIVERSITY OF OXFORD
1 SOUTH PARKS ROAD
OXFORD OX1 3TG
UNITED KINGDOM
E-MAIL: meinshausen@stats.ox.ac.uk

G. ROCHA
B. YU
DEPARTMENT OF STATISTICS
367 EVANS HALL
UNIVERSITY OF CALIFORNIA
BERKELEY, CALIFORNIA 94720-3860
USA
E-MAIL: gvrocha@stat.berkeley.edu
binyu@stat.berkeley.edu